\documentclass{amsart}
\usepackage{amsaddr}
\usepackage[english]{babel}
\usepackage{tablefootnote}
\usepackage{fancyhdr}

\usepackage{amsmath}
\usepackage{graphicx}

\usepackage{longtable}

\usepackage{amsmath, amsthm, amssymb}
\usepackage{color}
\usepackage{xspace}

\usepackage{enumitem}

\newtheorem{theorem}{Theorem}
\newtheorem{lemma}{Lemma}
\newtheorem{example}{Example}
\newtheorem{corollary}{Corollary}

\renewcommand{\>}{\rangle}

\theoremstyle{definition}

\newtheorem{definition}{Definition}
\newtheorem{remark}{Remark}

\def\sbsbs{\hsize=.33\textwidth\parindent=0pt\centering}
\long\def\sidebysidebyside#1#2#3{\figurecaptionfont
\hbox to\textwidth{\vtop{\sbsbs#1\vskip1sp}\hfill\vtop{\sbsbs#2\vskip1sp}\hfill\vtop{\sbsbs#3}}}

\definecolor{dkgreen}{rgb}{0,0.6,0}
\definecolor{gray}{rgb}{0.5,0.5,0.5}
\definecolor{mauve}{rgb}{0.58,0,0.82}

\definecolor{myblue}{rgb}{0.2,0.2,0.8}
\definecolor{mygreen}{rgb}{0.2,0.8,0.2}
\definecolor{myred}{rgb}{0.8,0.2,0.2}
\definecolor{mygold}{rgb}{0.6,0.4,0.2}
\definecolor{mypurple}{rgb}{0.6,0.2,0.4}
\definecolor{myteal}{rgb}{0.2,0.6,0.4}

\def\NN{\mathbb{N}}
\def\ZZ{\mathbb{Z}}

\newcommand{\mI}{\mathsf{I}}

\newcommand{\mL}{\mathsf{L}}
\newcommand{\mM}{\mathsf{M}}

\newcommand{\mO}{\mathsf{O}}

\newcommand{\mR}{\mathsf{R}}

\DeclareMathAlphabet\mathbfcal{OMS}{cmsy}{b}{n}

\usepackage[colorlinks=true, allcolors=blue]{hyperref}

\makeatletter
\patchcmd{\proof}{\let\@addpunct\@empty}{\let\@addpunct\@footnotemark\footnotetext{\@footnotetext}}{}{}
\makeatother

\makeatletter
\renewenvironment{theorem}[1][]
  {\ifstrempty{#1}{\trivlist\item[\hskip\labelsep{\bfseries Theorem.}]}{\trivlist\item[\hskip\labelsep{\bfseries Theorem (#1).}]}\itshape\footnotetext{\@footnotetext}}
  {\endtrivlist\@afterindentfalse}
\makeatother

\title{Rigid-Recurrent Sequences for Actions of Finite Exponent Groups}

\author{Cash Cherry}
\address{Department of Applied Mathematics and Statistics,\\ Colorado School of Mines,\\ Golden, Colorado 80401}

\begin{document}

\maketitle

\begin{abstract}
The focus of this paper is to better understand the coexistence of rigidity, weak mixing, and recurrence by constructing thin sets in the product of countably many copies of the finite cyclic group of order q. A Kronecker-type set K is a subset of this group on which every continuous function into the complex unit circle equals the restriction, to K, of a character in the group's Pontryagin dual. Ackelsberg proves that if, for all q $>$ 1, there exists a perfect Kronecker-type set generating a dense subgroup, then there exist large rigidity sequences for weak mixing systems of actions by countable discrete abelian groups.  Ackelsberg shows the existence of such sets for prime values of q, while we construct them for all q $>$ 1.
\end{abstract}

\let\thefootnote\relax
\let\thefootnote\relax

\addtocounter{footnote}{-2}

\section{Introduction}


$D_q$ is the countable product $\bigotimes_{n \in \NN} \ZZ/q\ZZ,$ equipped with the product topology. $K_q$ sets are subsets of $D_q$ on which every continuous function into the unit circle equals the restriction of a character in $D_q$'s Pontryagin dual $\widehat{D}_q$. We construct perfect sets of type \emph{$K_q$} which generate dense subgroups of $D_q$ to demonstrate that weak mixing actions by countable discrete abelian groups of finite exponent can exhibit recurrence.  Here we repeat definitions from Ackelsberg \cite{Ackelsberg}:

\begin{definition}
    Let $\Gamma$ be a countable discrete abelian group. A \emph{measure preserving system} is a quadruple $(X,\mathcal{B},\mu,(T_g)_{g \in \Gamma}),$ where $(X,\mathcal{B},\mu)$ is a non-atomic Lebesgue probability space, and $(T_g)_{g \in \Gamma}$ is an action of $\Gamma$ by measure-preserving transformations. A sequence $(a_n)_{n \in \NN} \subseteq \Gamma$ is \emph{rigid} for the system if for every $f \in L^2(\mu),$ $\| f\circ T_{a_n} - f \|_2 \to 0.$  
\end{definition}

\begin{definition}
    A sequence $(\Phi_N)_{N \in \NN}$ of finite subsets of $\Gamma$ is a \emph{Følner sequence} if for every $x \in \Gamma,$
    $$
    \frac{|(\Phi_N + x) \triangle \Phi_N|}{|\Phi_N|} \xrightarrow[N \to \infty]{} 0
    $$
    A system $(X,\mathcal{B},\mu,(T_g)_{g \in \Gamma})$ is \emph{weak mixing} if for every Følner sequence $(\Phi_N)_{N \in \NN}$ in $\Gamma$ and all $A,B \in \mathcal{B},$
    $$
    \frac{1}{|\Phi_N|} \sum_{g \in \Phi_N} |\mu(A \cap T_g B) - \mu(A)\mu(B)| \xrightarrow[N \to \infty]{} 0
    $$
\end{definition}

\begin{definition}
    Let $\Gamma$ be a countable discrete abelian group. A set $R \subseteq \Gamma$ is a \emph{set of recurrence} if for every measure-preserving system $(X,\mathcal{B},\mu,(T_g)_{g \in \Gamma})$ and every $A \in \mathcal{B}$ with $\mu(A)>0,$ there exists $r \in R\setminus \{0\}$ such that $\mu(A \cap T_r^{-1} A) > 0.$
\end{definition}

\begin{definition}
    Let $\Gamma$ be a countable discrete abelian group. A sequence $(r_n)_{n \in \NN}$ in $\Gamma$ is \emph{rigid-recurrent} if $(r_n)_{n \in \NN}$ is rigid for some weak mixing measure preserving system and $\{ r_n : n \in \NN \}$ is a set of recurrence. Such a sequence is furthermore said to be \emph{freely rigid-recurrent} if it is rigid for a \emph{free} measure-preserving system, that is, a system $(X,\mathcal{B},\mu,(T_g)_{g \in \Gamma})$ for which
    $$\mu(\{ x \in X : T_g x = x \}) = 0$$
    when $g \neq 0.$
\end{definition}

Ackelsberg conjectures (Conjecture 1.6 in \cite{Ackelsberg}) that every countable discrete abelian group $\Gamma$ contains a sequence $S$ and a finite index subgroup $\Delta \leq \Gamma$ such that every translate of $S$ by an element of $\Delta$ is a freely rigid-recurrent sequence. Griesmer \cite{griesmer} had previously shown this for the case where $\Gamma = \Delta = \ZZ.$ Section 7.2 of \cite{Ackelsberg} provides a sufficient condition for this proposition in the existence of perfect Kronecker or $K_q$ sets which generate dense subgroups of locally compact abelian groups. There, the conjecture 
is proven
for the case where $\Delta$ is of the form $\bigoplus_{j=1}^N \bigoplus_{n=1}^\infty \ZZ/{p_j}\ZZ$ for distinct primes $p_1,p_2, \dots p_N$, and $\Gamma = \Delta \bigoplus \mathrm{F}$ for some finite abelian group $\mathrm{F}$. Using Ackelsberg's characterization, we extend this to all countable discrete abelian groups of finite exponent.

\section{Background}

Fixing $q \in \NN,$ we define the following notation:

\begin{longtable}{l p{300pt}}
    $\NN_n$ & The finite subset of $\NN$ given by $\{1,2,\dots,n\}.$ \\
    $C_q$ & The ring $\ZZ/q\ZZ$ with operations given by addition and multiplication mod $q.$ \\
    $C_q^n$ & The $C_q$ module given by functions $\NN_n \to C_q,$ represented by column vectors with entries in $C_q.$\\
    $C_q^{n \times m}$ & The set of $n\times m$ matrices with entries in $C_q.$\\
    $D_q$ & The $C_q$ module given by functions $\NN \to C_q.$ The additive group is isomorphic to the standard definition $D_q = \bigotimes_{n \in \NN} \ZZ / q \ZZ .$\\ 
    $x \mapsto x_{|n}$ & The homomorphism $D_q \to C_q^n$ given by restricting the domain of $x \in D_q$ from $\NN$ to $\NN_n.$  
    \\
    $K_{|n}$ & Given $K \subset D_q,$ the set $\{x_{|n} : x \in K\}$
    \\
    $[t]$ & Given $t \in C_q^n,$ this is the \emph{cylinder set} $\{x \in D_q : x_{|n} = t\}.$ We give $D_q$ the product topology with cylinder sets as basic open sets.\\
    $\kappa_f$ & For a function $f: D_q \to C_q,$ this is the \emph{continuity index}, defined to be $\min \{ \kappa \in \NN : x_{|\kappa} = y_{|\kappa} \hspace{5pt}\text{implies} \hspace{5pt}f(x) = f(y)  \}.$ $f$ is continuous under the product topology exactly when such a continuity index exists.
    \\
    $\widehat{D}_q$ & This is the group (under addition) of continuous homomorphisms from the additive group of $D_q$ into the additive group of $C_q.$ Elements of $\widehat{D}_q$ are called \emph{characters.}\\
    $e_i$ & This is the element of $C_q^n$ or $D_q$ defined by $e_i(j) = \begin{cases}
    1 & \text{if}\hspace{5pt} i = j\\
    0 & \text{if}\hspace{5pt} i \neq j\end{cases}.$\\
    $\mO$ & A matrix or vector of all zeros in $C_q$ whose size is determined by context.
\end{longtable}

\begin{remark} 
$\widehat{D}_q$ as defined above is isomorphic to the Pontryagin dual group of the additive group of $D_q.$ The difference amounts to a choice between multiplicative and additive notation. Additive notation is chosen for this particular construction because of connections between the desired properties and $C_q-$valued matrices.
\end{remark}

\begin{definition}
    A set $K \subset D_q$ is of \emph{type} $K_q$ if, for all continuous $f : K \to C_q,$ there exists a character $\chi \in \widehat{D}_q$ such that $\chi_{|K} = f.$ These are also referred to as \emph{$K_q$ sets}, or sets of \emph{Kronecker-type}.
\end{definition}

\begin{definition}
    A set $K \subset D_q$ \emph{topologically generates} $D_q$ if $K$ algebraically generates a dense subgroup of $D_q.$
\end{definition}

\begin{definition}
    For a given $q,$ $D_q$ is \emph{topologically Kronecker-generated} if there exists a perfect $K_q$ set $K \subset D_q$ such that $K$ topologically generates $D_q.$
\end{definition}    

\section{Results}

\subsection{Characterization of Perfect Topological Generating sets of Type $K_q$ in $D_q$} \hfill

Consider a perfect set $K \subset D_q.$ For any $n \in \NN,$ let $k_{n} = |K_{|n}|,$ and choose  $\mM_n \in C_q^{n \times k_n}$ to be any matrix whose set of columns equals $K_{|n}.$ 
We characterize sets of type $K_q$ and sets which topologically generate $D_q$ in terms of the left and right invertibility of the corresponding matrices $\mM_n.$ 

\begin{remark}
    As left and right invertibility are invariant under permutations of columns, this arbitrary choice of $\mM_n$ makes no difference to the following results.
\end{remark}

\begin{lemma}\label{charvec}  
    Each character $\chi \in \widehat{D}_q$ has a vector representation $c \in C_q^{\kappa_\chi}$ satisfying $x_{|\kappa_\chi}^T c = \chi(x)$ for all $x \in D_q.$ It is given by letting $c_i = \chi(e_i).$

\begin{proof}
    Given any $x \in D_q,$
    $$
    x_{|\kappa_{\chi}}^T c = \sum_{i=1}^{\kappa_{\chi}} x_i \chi(e_i)  = \chi\Bigg{(} \sum_{i=1}^{\kappa_{\chi}} x_i e_i \Bigg{)} = \chi(x_{\kappa_{\chi}}) = \chi(x)
    $$
\end{proof}
\end{lemma}

\begin{lemma}\label{K_q}
    A perfect set $K \subset D_q$ is $K_q$ if $\mM_n$ has a left inverse $\mL_n$ for infinitely many $n.$
\begin{proof}
    Suppose $\mM_n$ is left invertible for infinitely many $n.$
    Let $f: K \to C_q$ be continuous.
    Fix $n\geq \kappa_f$ such that $\mM_n$ has left a inverse $\mL_n$.  Create the column vector $\bar{f}(\mM_n^T) \in C^{k_n}$ so it has entries given by applying $f$ to elements of $K$ which truncate to the corresponding rows of $\mM_n^T.$ That is, picking some $x^{(1)},\dots,x^{(k_n)} \in K$ such that
    $$\mM_n^T = \begin{bmatrix}
        {x^{(1)}_{|n}}^T\\ \vdots \\{x^{(k_n)}_{|n}}^T
    \end{bmatrix}, \hspace{10pt} \text{we can write} \hspace{10pt} \bar{f}(\mM_n^T) = \begin{bmatrix}
        f(x^{(1)})\\ \vdots \\f(x^{(n)})
    \end{bmatrix}$$
    Consider the character $\chi$ having 
    $$ \chi(x) = (\mL_n x_{|n})^T \bar{f}(\mM_n^T)$$
    As $\mL_n \mM_n = \mI_k,$  
    it must be the case that $\mL_n x_{|n}$ is a standard basis vector corresponding to the row of $\mM_n^T$ containing $x_{|n}.$ Thus, by the definition of $\bar{f}(\mM_n^T),$ it follows that $\chi(x) = (\mL_n x_{|n})^T \bar{f}(\mM_n^T) = f(x)$ for all $x \in K.$ As we have arbitrarily chosen a continuous $f: K \to C_q,$ and constructed a character $\chi \in \widehat{D_q}$ satisfying $\chi_{|K} = f,$ it follows that $K$ is $K_q.$
\end{proof}
\end{lemma}

\begin{remark}
    A statement similar to the converse of Lemma \ref{K_q} is true: If $K$ is $K_q,$ then for all continuous $f: K \to C_q,$ there exists $n \in \NN$ such that $\mM_n^T c = \bar{f}[\mM_n^T]$ has a solution $c \in C_q^n.$ Though this looser condition is implied by left invertibility of $\mM_n$ for infinitely many $n,$ the converse of Lemma \ref{K_q} does not hold, as can be seen in Lemma \ref{1}.
\end{remark}

\begin{example} \label{1}
    Letting $K = \{e_n : n \in \NN\},$ it is clear that $K$ is of type $K_q$: Given a continuous $f : K \to C_q,$ we can define a character $\chi$ which sends \newline $x \mapsto x_{\kappa_f}^T [f(e_1),\dots,f(e_{\kappa_f})]^T,$ and it will be the case that $\chi_{K} = f.$ However, every $\mM_n$ has $n$ columns for each $e_i$ with $i \leq n,$ and a column of zeros for each $e_i$ with $i > n.$ Thus it will be the case that $\mM_n \in C_q^{n \times n + 1}$ for all $n,$ and thus, by its dimension cannot be left invertible. Therefore, $\mM_n$ is left invertible for no $n,$ but $K$ is of type $K_q.$ This is a counterexample to the converse of Lemma \ref{K_q}.
\end{example}

\begin{lemma}\label{generates}
    A perfect set $K \subseteq D_q$ topologically generates $D_q$ if and only if 
    $\mM_n$ has a right inverse for infinitely many $n.$
\begin{proof}
    Suppose $\mM_n$ is right invertible for infinitely many $n.$ 
    Pick $m \in \NN,$ and a cylinder set $[t]$ with $t \in C_q^m.$ Choosing $n\geq m$ such that $\mM_n$ has right inverse $\mR_n$, we can let $$\alpha = \mR_n \begin{bmatrix}t\\\mO\end{bmatrix}, \hspace{10pt} \text{so that} \hspace{10pt} \mM_n \alpha = \begin{bmatrix}t\\\mO\end{bmatrix}$$
    Picking $x^{(1)},\dots,x^{(k_n)} \in K$ as in the proof of Lemma \ref{K_q}, let $y = \sum_{i=1}^{k_n} \alpha_i x^{(i)}.$ It follows from construction that $y \in \<K\>$ and $y \in [t].$ As $[t]$ was chosen arbitrary, this shows that $\<K\>$ is dense in $D_q,$ that is, $K$ topologically generates $D_q.$

    Now, suppose that $K$ topologically generates $D_q.$ Picking $n \in \NN,$ we can find a solution $r_i$ to
    $\mM_n r_i = e_i$
    for all $e_i \in C_q^{k_n},$ $1 \leq i \leq n.$  Otherwise, there would be some $e_i \in C_q^{k_n}$ with $[e_i] \cap \<K\> = \emptyset,$ and $K$ would not topologically generate $D_q$ (Note that for this reason, it must be the case that $n \leq k_n$). Thus, the equation
    $$
    \mM_n \mR_n = \mI
    $$
    has a solution for $\mR_n \in C_q^{k_n \times n}$ whose $i$th column is given by $r_i.$
\end{proof}
\end{lemma}

\subsection{Topological Kronecker Generation of \texorpdfstring{$D_q$}{Dq}}

\begin{theorem}
    $D_q$ is topologically Kronecker generated for all $q \in \NN \setminus\{1\}.$
\begin{proof}
    We define the sequence of $C_q-$valued matrices $2^n \times 2^n$ matrices $\mM_{2^n}$ by letting $\mM_{1} = [1]$ and  defining
    $$
    \mM_{2^n} = \begin{bmatrix}
    \mM_{2^{n-1}}&\mM_{2^{n-1}}\\
    \mO&\mI
    \end{bmatrix}
    $$
    for all $n \in \NN.$ The first few such matrices are $\mM_{1} = [1]$,
    $$
    \mM_{2} = \begin{bmatrix}
        1&1\\0&1
    \end{bmatrix}, 
    \mM_{4} = \begin{bmatrix}
        1&1&1&1\\
        0&1&0&1\\
        0&0&1&0\\
        0&0&0&1
    \end{bmatrix},
    \mM_8 = \begin{bmatrix}
        1&1&1&1&1&1&1&1\\
        0&1&0&1&0&1&0&1\\
        0&0&1&0&0&0&1&0\\
        0&0&0&1&0&0&0&1\\
        0&0&0&0&1&0&0&0\\
        0&0&0&0&0&1&0&0\\
        0&0&0&0&0&0&1&0\\
        0&0&0&0&0&0&0&1
    \end{bmatrix}
    \dots
    $$
    Note that each $\mM_{2^n}$ is upper unitriangular and thus has determinant 1. It follows that each $\mM_{2^n}$ invertible over $C_q.$  
    
    \noindent Letting $S_{2^n}$ be the set of columns of $\mM_{2^n},$ define
    $$K := \bigcap_{n \in \NN_0} \bigcup_{t \in S_{2^n}} [t]$$
    Examining $\mM_8$ shows us that the element of finite support $e_1 + e_2 + e_6$ is contained in $K.$ However, note that $K$ also contains elements of infinite support such as $\sum_{n=0}^\infty e_{2^n}.$

    This is a perfect set: By the recursive construction of $\mM_{2^n},$ the columns of $\mM_{2^n}$ are all columns of $\mM_{2^{n-1}}$ with $2^{n-1}$ zeros appended below, and all columns of $\mM_{2^{n-1}}$ with some column
    from $\mI_{2^{n-1}}$ appended below. This is a binary choice at each $n$ leading to distinct cylinder sets, allowing us to place the cylinder sets on the nodes of an infinite complete binary tree which is disjoint at each level, and obeys an inclusion ordering for which $K$ is the intersection of the union of each level. Since the cylinder sets are uncountable, compact, and contain no isolated points, it follows that $K$ is a perfect set. 

    Furthermore, $\mM_{2^n}$ as defined here exactly satisfies the definition in Lemma \ref{K_q}. Thus, we can use the invertibility of $\mM_{2^n}$ for all $n$ to apply Lemmas \ref{K_q} and \ref{generates} and conclude that $K$ is also a  $K_q$ set which topologically generates $D_q.$ It follows by definition that $D_q$ is topologically Kronecker generated.    
\end{proof}
\end{theorem}

From Proposition 7.10 and Theorem 7.11 in Ackelsberg \cite{Ackelsberg}, we get the following corollary:

\begin{corollary}
    For every countable discrete abelian group $\Gamma$ with finite exponent, there exists a sequence $(r_n)_{n \in \NN}$ in $\Gamma$ and a finite index subgroup $\Delta \leq \Gamma$ such that for every $s \in \Delta,$ $(r_n - s)_{n \in \NN}$ is freely rigid-recurrent sequence.
\end{corollary}

\noindent \textbf{Acknowledgements.} I am immensely grateful to John Griesmer for introducing me to these topics.

\bibliographystyle{alpha}


{{\footnotesize 
\textit{E-mail address:} \href{ccherry@mines.edu}{ccherry@mines.edu}
}}
\end{document}